\documentclass[12pt,a4paper]{article}
\usepackage[ngerman]{babel}
\usepackage[margin=4cm]{geometry}
\usepackage[utf8]{inputenc}
\usepackage{avant}
\usepackage{graphicx}
\usepackage{pslatex}
\usepackage{textcomp}
\usepackage{gensymb}
\usepackage{setspace}
\usepackage{lmodern}
\usepackage{amsmath,amssymb, amsthm}
\usepackage{cases}

\newlength{\bibitemsep}
\newlength{\bibparskip}
\let\oldthebibliography\thebibliography
\renewcommand\thebibliography[1]{%
  \oldthebibliography{#1}%
  \setlength{\parskip}{\bibitemsep}%
  \setlength{\itemsep}{\bibparskip}%
}

\numberwithin{equation}{section}
\newtheorem{theorem}{Theorem} [section]
\newtheorem{lemma}[theorem]{Lemma} 
 
\newtheorem{cor}[theorem]{Corollary}
\begin{document}
\pagenumbering{arabic}

\title{Moment Generating Stirling Numbers and Applications}
\author{Ludwig Frank 
\\ {\small Faculty of Computer Science, Technical University of Applied Sciences,} 
\\ {\small Hochschulstraße 1, D-83024 Rosenheim, Germany,}
\\ {\small ludwig.frank@th-rosenheim.de}
}
\date{}
\maketitle
\begin{quote}
\begin{small}
\textbf{Abstract.} 
In this paper, we investigate certain combinatorial numbers, the \textit{moment generating Stirling numbers}. They are a special case of Hsu's generalized Stirling numbers and satisfy many more properties and combinatorial identities than are known in the general case. As application, we provide the computation of the moments and central moments of the phase type distribution, the recurrence time in Markov chains, the geometric distribution, the negative binomial distribution and of a class of distributions generalizing the negative binomial distribution. All computations can be performed in closed form without recursion.  We also present the relationship to the Markov renewal process.   
\end{small}



\textbf{Statements and Declarations.} No funds, grants, or other support was received.

\end{quote}

\section{Introduction}\  

Stirling numbers and their extensions play an important role in combinatorics, probability theory and graph theory (cf. e.g. \cite{Hertz}). This has lead to many generalizations. A comprehensive summary of them can be found in \cite{Hsu}, where Stirling numbers of the first and second kind, the Lah numbers and numerous variants of them were put together into a single parametrized formula. Despite the generality of this approach, new extensions of the Stirling numbers have been found, cf. e.g. \cite{Belbachir}, \cite{Belkhier}, \cite{Cakic}, \cite{Choi}, \cite{Lang}, \cite{Mansour}, \cite{Mohr} or \cite{Steingrimsson}. In this paper a special case of \cite{Hsu} is examined in detail. It turns out that it satifies many more properties and identities than are known in the general case, and that these are useful in applications. 
We call them \textit{moment generating Stirling numbers} (\textit{MSN}), because the easiest way to show their benefits is to compute the moments of some probability distributions by means of closed formulae and not by the successive differentiation of a generating function. 

In Section 2.1 we present the formal definition, a combinatorial interpretation and the relation to the \textit{D-Stirling numbers} and the \textit{degenerate weighted} Stirling numbers, which are also a special case of the concept in \cite{Hsu}. Then, the fundamental results are presented, from which all further properties can be derived, and some often used special values. There are two types of formulas, those which generalize similar formulas of the corresponding ones of the Stirling numbers of the second kind, and those which have no counterpart, because they involve the additional new parameter. The structure of the subsection arises from the fact that we summarize the results according to the summation parameter. Finally, we investigate the existence of inverses. In Section 2.2 we show the representation of the \textit{MSNs} by multinomial coefficients, and we present various generating functions that make it possible to find further identities between the \textit{MSNs}. In Section 3, we show how \textit{MSNs} can be applied. We present formulas for the moments of the phase type distribution, the recurrence time in Markov chains, the geometric, the negative binomial distribution, a class of distributions generalizing the negative binomial distribution and the relationship to the Markov renewal process. Also the application on the central moments is shown in Section 3. Finally, we list some open questions for further investigations.
\section{Moment Generating Stirling Numbers of the second kind}
\subsection{Definitions and Basic Concepts} 

For integers $i \ge 0$, $j \ge 0$ and real $k$, we define the \textbf{\textit{moment generating Stirling numbers (of the second kind), "`MSN"'}},
\begin{equation}
	b_{i,j,k}=\sum_{r=0}^{j} \binom{j}{r} \cdot (-1)^{j-r} \cdot (r+k)^i.
	\label{a1}
\end{equation}

It is noticeable that we do not divide by the factorial of the second parameter as in most of the other variants, but the formulas for the identities and the applications become simpler this way.

For integer $k \ge 0$, there is a combinatorial interpretation of $b_{i,j,k}$: It is the number of ways to put $i$ distinguishable balls into $j+k$ distinguishable boxes so that none of the boxes 1 to $j$ are left empty: This event is $\overline{\cup_{r=1}^{j} E_r}$, where $E_r$ denotes the event that box $r$ remains empty. The sieve formula yields 
\begin{align*}
	|\overline{\cup_{r=1}^{j} E_r}|=\sum_{r=0}^{j} \binom{j}{r} \cdot(-1)^r \cdot (k+j-r)^i=b_{i,j,k},
\end{align*}
because $|E_{i_1} \cap \ldots \cap E_{i_r}|=(k+j-r)^i$. After Theorem \ref{Interpretation} a second interpretation is possible.

The \textit{D-Stirling numbers} in \cite{Huang} allow for each box individual lower and also upper bounds for the number of balls. For \textit{MSNs} there are no upper bounds and the lower bounds are 1 and 0. For (ordinary) Stirling numbers the lower bound is always 1. In Hsu's and Shiue's generalized Stirling numbers the \textit{MSNs} (divided by $j!$) correspond to the parameter triple $(\alpha, \beta, r)=(0,1,k)$. For the \textit{degenerate weighted} Stirling numbers (cf. e.g. \cite{Cakic} or \cite{Hsu}), $S(n,k,\lambda|\Theta)=b_{n,k,\lambda}/k!$ with $\Theta=0$ holds. 

{\allowdisplaybreaks
The following lemma summarizes some obvious and often used properties:
\begin{lemma} \label{lemma1}
For integers $i \ge 0$, $j \ge 0$ and real $k$, the following properties hold:
\begin{align}
  \text{a) }  &b_{0,j,k}=0, \text{ for } j > 0, \label{a2}\\
  \text{b) }  &b_{0,0,k}=1,  \label{a3}\\
	\text{c) }  &b_{i,0,k}=k^i, \label{a4}\\
	\text{d) }  &b_{i,j,0}=S_{i,j} \cdot j! 
	\text{ with the Stirling numbers }S_{i,j} \text{ of the second kind,}\label{a5}\\
	\text{e) }  &b_{i,1,k}=(k+1)^{i}-k^i.\label{a6}
\end{align}
\end{lemma}
The proofs follow immediately from the definition in (\ref{a1}) and the definition of the Stirling numbers (cf. e.g. \cite{Heise}). 
Item d) allows to translate all known results about Stirling numbers of the second type to the formulae for $b_{i,j,0}$ and vice versa. Further, d) has the consequence that $b_{i,j,0}$ is the number of all surjective functions of an $i$-set onto a $j$-set (cf. e.g. \cite{Heise}). 

After the next theorem, which contains the basic important results about \textit{MSNs}, further simple results like in Lemma \ref{lemma1} can be proved. 
\begin {theorem} \ 
For integers $i \ge 0$, $j \ge 0$ and real $k$, the following properties hold:
\begin{align}
&\text{a)  } b_{i,j,k+1}=b_{i,j,k}+b_{i,j+1,k}, \label{a8}\\
&\text{b)  } b_{i+1,j+1,k}=(j+1) \cdot b_{i,j,k} + (j+1+k) \cdot b_{i,j+1,k},\label{a9} \\
&\text{c)  } b_{i+1,j+1,k}=(j+1) \cdot b_{i,j,k+1} + k\cdot b_{i,j+1,k}. \label{a12} 
\end{align}
\end{theorem}
\textit{Proof.} a) follows from the definition in (\ref{a1}).
b) Due to (\ref{a1}),
\begin{align*}
&(j+1) \cdot b_{i,j,k} + (j+1+k) \cdot b_{i,j+1,k} = \\
&\sum_{r=0}^{j} (-1)^{j-r} \cdot \left((j+1) \cdot \binom{j}{r} -(j+1+k) \cdot \binom{j+1}{r}\right) \cdot(r+k)^i + \\
&	\ \ \ 	(j+1+k) \cdot (j+1+k)^i =\\
&\sum_{r=0}^{j+1} (-1)^{j+1-r} \cdot \binom{j+1}{r} \cdot (r+k)^{i+1}. 
\end{align*}
c) follows from (\ref{a9}) and (\ref{a8}).$\qed$

Since $S_{i,j} \ge 0$, (\ref{a5}) and (\ref{a8}) show that $b_{i,j,k} \ge 0$ for $k \ge 0$, but for negative values of $k$ $b_{i,j,k}$ can also be negative, e.g. $b_{1,0,-1}=-1$.

The following corollary contains the simple properties announced after Lemma \ref{lemma1}.
\begin{cor} For integers $i \ge 0$, $j \ge 0$ and real $k$ 
\begin{align}
&\text{ a) } \text{If } i<j, \text{ then } b_{i,j,k}=0, \label{a13} \\
&\text{ b) } b_{i,i,k}=i!, \label{a14} \\
&\text{ c) } b_{i+1,i,k}=(i+1)! \cdot \frac{i+2 \cdot k}{2}, \label{a15} \\
&\text{ d) } (j+1) \cdot b_{i,j,1}=b_{i+1,j+1,0}, \label{a16} \\
&\text{ e) } b_{1,j,k}=
\begin{cases}
    1 & \text{, if } j=1\\
		k & \text{, if } j=0\\
		0 & \text{, if } j>1,
\end{cases} \label{a7} \\
&\text{ f) } b_{i,j,k-j}=(-1)^{i+j} \cdot b_{i,j,-k }, \label{n1} \\
&\text{ g) } b_{i,2 \cdot k,-k}= 0, \text{ if } k \ge 0 \text{ and } i \text{ is odd}. \label{n2}
\end{align}
\end {cor}
\textit{Proof.} Due to (\ref{a2}), (\ref{a3}) resp. (\ref{a4}), the proofs of a) until c) can be performed by means of induction on $i$. 
d) follows from (\ref{a12}).
e) follows from (\ref{a6}) for $j=1$, from (\ref{a4}) for $j=0$ and from (\ref{a13}) for $j>1$.
f) follows from the definition in (\ref{a1}) by means of index transformation and 
g) immediately from f). $\square{}$ 

The following corollary starts the series of theorems and corollaries, in which over the  first parameter is summed. 

\begin{cor} For integers $i \ge 0$, $j \ge 0$ and real $k$
\begin{align}
& \text{a)} &b_{i,j+1,k}=&\sum_{r=0}^{i-1} \binom{i}{r} \cdot b_{r,j,k}, \label{a10}\\
& \text{b)} &b_{i,j,k+1}=&\sum_{r=0}^{i} \binom{i}{r} \cdot b_{r,j,k},\label{a11}\\ 
& \text{c)} &b_{i+1,j+1,k}=&(j+1)\cdot \sum_{r=0}^{i} \binom{i}{r} \cdot b_{r,j,k}+k  
\cdot b_{i,j+1,k}, \label{a12a} \\
& \text{d)} &b_{i+1,j+1,k}=&(j+1) \cdot \sum_{r=j}^{i} b_{r,j,k} \cdot (j+k+1)^{i-r},\label{a12b} \\
& \text{e)} &b_{i+1,j+1,k}=&(j+1) \cdot \sum_{r=j}^{i} b_{r,j,k+1} \cdot k^{i-r}.\label{n3}
\end{align}
\end{cor}
\textit{Proof.} a) Due to the definition
\begin{align*}
		&b_{i,j+1,k}=\sum_{r=0}^{j+1} \binom{j+1}{r} \cdot (-1)^{j+1-r} \cdot (r+k)^i= \\
		&\sum_{r=0}^{j+1} \binom{j}{r} \cdot (-1)^{j+1-r} \cdot (r+k)^i+  \sum_{r=1}^{j+1} \binom{j}{r-1} \cdot (-1)^{j+1-r} \cdot (r+k)^i= \\
		&-b_{i,j,k}+\sum_{s=0}^{i} \binom{i}{s} \cdot \sum_{r=0}^{j}\binom{j}{r} \cdot (-1)^{j-r} \cdot (r+k)^s= 
		\sum_{s=0}^{i-1} \binom{i}{s} \cdot b_{s,j,k}.
\end{align*}
b) is a consequence of (\ref{a10}) and (\ref{a8}).
c) follows immediately from (\ref{a12}) and (\ref{a11}). 
d) Induction on $i$ by means of (\ref{a7}), (\ref{a3}) and (\ref{a9}).
e) Induction on $i$ by means of (\ref{a9}).$\square{}$

Item d) is known from \cite{DLMF}, (formula 26.4.24), if $k=0$ and $b_{i,j,0}$ is replaced by $S_{i,j} \cdot j!$ (cf. (\ref{a5})). Like (\ref{a8}), Item b) can help to compute the \textit{MSNs} successively for growing $k \ge0$, while Item e), on the other hand, is suitable for computing the \textit{MSNs} for decreasing $k <0$, starting with $b_{i,j,0}$.

The following theorem, also known from \cite{DLMF} for $k_1=k_2=0$, and its corollary continue the summation over the first parameter of the \textit{MSNs} in a form which reminds of convolution theorems.
\begin {theorem}\ 
For integers $i \ge 0$, $j_1 \ge 0$, $j_2 \ge 0$ and real $k_1$ and $k_2$
\begin{equation}
b_{i,j_1+j_2,k_1+k_2}=\sum_{r=0}^{i} \binom{i}{r} \cdot b_{r,j_1,k_1}\cdot b_{i-r,j_2,k_2},\label{a17}\\ 
\end{equation}
\end{theorem}
\textit{Proof.} Due to the definition in (\ref{a1})
\begin{align*}
&\sum_{r=0}^{i} \binom{i}{r} \cdot b_{r,j_1,k_1}\cdot b_{i-r,j_2,k_2}= \\
&\sum_{r=0}^{i} \binom{i}{r} \cdot \sum_{s=0}^{j_1}  \binom{j_1}{s} \cdot (-1)^{j_1-s} \cdot (s+k_1)^{r} \cdot \sum_{s=0}^{j_2}  \binom{j_2}{s} \cdot (-1)^{j_2-s} \cdot (s+k_2)^{i-r}=\\
&\sum_{l=0}^{j_1+j_2} \sum_{s=0}^{l} \binom{j_1}{s} 
 \cdot \binom{j_2}{l-s}\cdot (-1)^{j_1+j_2-l }\cdot \sum_{r=0}^{i} \binom{i}{r} \cdot (s+k_1)^r \cdot (l-s+k_2)^{i-r}=\\
&\sum_{l=0}^{j_1+j_2} \sum_{s=0}^{l} \binom{j_1}{s} 
 \cdot \binom{j_2}{l-s}\cdot (-1)^{j_1+j_2-l }\cdot (l+k_1+k_2)^{i}=b_{i,j_1+j_2,k_1+k_2}. \ \square{}
\end{align*}

\begin{cor} For integers $i \ge 0$, $j \ge 0$, $j_1 \ge 0$, $j_2 \ge 0$ and real $k_1$, $k_2$ and $k$ 
\begin{align}
&  \text{a) }  &b_{i,j_1+j_2,0}=&\sum_{r=0}^{i} \binom{i}{r} \cdot b_{r,j_1,k}\cdot b_{i-r,j_2,-k},\label{a18}\\
&  \text{b) }  &b_{i,j,0}=&\sum_{r=0}^{i} \binom{i}{r} \cdot k^r \cdot b_{i-r,j,-k},\label{a19}\\  
&	\text{c) }  &b_{i,j,k}=&\sum_{r=0}^{i} \binom{i}{r} \cdot k^{i-r} \cdot b_{r,j,0},\label{a20}\\  
&	\text{d) }  &b_{i,j,k_1+k_2}=&\sum_{r=0}^{i} \binom{i}{r} \cdot k_1^{i-r} \cdot b_{r,j,k_2},\label{a21}\\  
&	\text{e) }  &j \cdot b_{i,j-1,k+1}=&\sum_{r=0}^{i} \binom{i}{r} \cdot k^{i-r}\cdot b_{r+1,j,0}.\label{a23} 
\end{align}
\end{cor}
\textit{Proof.} a) follows from (\ref{a17}) with $k_1=-k_2=k$. b) follows from (\ref{a18}) with $j=j_2$ and $j_1=0$, because (cf. (\ref{a4})) $b_{r,0,k}=k^r$. c) follows from (\ref{a17}) with $k_1=0$, $k_2=k$, $j_2=0$ and the fact that $b_{i-r,0,k_1}=k_1^{i-r}$. d) follows from (\ref{a17}) with $j_2=0$ for the same reason. e) follows from (\ref{a21}) with $k_2=1$ and (\ref{a16}). $\square{}$

\paragraph{}
By means of (\ref{a5}) the Items c) and e) show tight connections between the Stirling numbers of the second type and the \textit{MSNs} with arbitrary $k$ (cf. \cite{DLMF}). Furthermore, the connection between $b_{i,j,k}$ and $b_{i,j,-k}$ becomes visible. 

The series of theorems and corollaries with summation over the second parameter begins with the following theorem.  
\begin {theorem}\ 
For integers $i \ge 0$, $n \ge 0$ and real $k$
\begin{equation}
(n+k)^i= \sum_{r=0}^{i} \binom{n}{r} \cdot b_{i,r,k}. \label{a7a}
\end{equation}
\end{theorem}
\textit{Proof.} Due to the inversion formula for binomial coefficients (cf. e.g. \cite{Heise}) and the definition of the \textit{MSNs} in (\ref{a1}), we get for $n \ge 0$ 
	\[\sum_{r=0}^{n} \binom{n}{r} \cdot b_{i,r,k} = 
	\sum_{l=0}^{n} \sum_{r=l}^{n} \binom{n}{r} \cdot \binom{r}{l} \cdot (-1)^{r-l} \cdot 
	(l+k)^i= (n+k)^i,
\]
and due to $b_{i,r,k}=0$ for $r>i$ (cf. (\ref{a13})) the assertion (\ref{a7a}). \ $\square{}$ 
\paragraph{}
For $k=0$, the preceding theorem is known in the form
	\[n^i= \sum_{r=0}^{i} \binom{n}{r} \cdot b_{i,r,0} = 
	\sum_{r=0}^{i} \binom{n}{r} \cdot S_{i,r} \cdot r!= \sum_{r=0}^{i} S_{i,r} \cdot (n)_{r}\
\]
with the Stirling numbers $S_{i,r}$ of the second type and $(n)_{r}=n \cdot (n-1)  \cdot \ldots \cdot (n-r+1)$ (cf. e.g. \cite{Heise}). As in the case $k=0$, (\ref{a7a}) remains valid, if $n$ is replaced by any $x \in \mathbb{R}$. The following corollary follows immediately:

\begin{cor} For integers $i \ge 0$, $n \ge 0$, $l$ and real $k$ 
\begin{align}
  \text{a) }  &(n+k)^i= \sum_{r=0}^{i} \binom{n+l}{r} \cdot b_{i,r,k-l} , 
	\text{ if } n+l \ge 0 \label{a7a1}\\
	\text{b) }  &(l \cdot n+k)^i=\sum_{r=0}^{i} \binom{n}{r} \cdot b_{i,r,(l-1)\cdot n +k},\text{ if } l \ge 0. \label{a7a2}
\end{align}
\end{cor}
\textit{Proof.} a) follows from (\ref{a7a}) due to $n+k=(n+l)+(k-l)$ and b) by means of (\ref{a21}) due to 
\begin{align*}
	&(l \cdot n+k)^i=\sum_{s=0}^{i} \binom{i}{s} \cdot (n+ (l-1) \cdot n)^s \cdot k^{i-s} =\\
	&\sum_{s=0}^{i} \binom{i}{s} \cdot  k^{i-s} \sum_{r=0}^{s} \binom{n}{r} \cdot b_{s,r,(l-1)\cdot n}=
	\sum_{r=0}^{i} \binom{n}{r} \cdot b_{i,r,(l-1)\cdot n+ k}. \ \square{}
	\end{align*}
	
\begin {theorem}\ 
For integers $i \ge 0$, $j \ge 0$ and real $k$
\begin{equation}
\sum_{r=0}^{j} (-1)^{j-r} \cdot b_{i,r,k}=b_{i,j+1,k-1}+ (-1)^{j} \cdot (k-1)^i.\label{a24} 
\end{equation}
\end{theorem}
\textit{Proof.} Due to (\ref{a21}) with $k=k_1$ and $k_2=0$ and the definition in (\ref{a1})
\begin{align*}
&\sum_{r=0}^{j} (-1)^{j-r} \cdot b_{i,r,k}=
\sum_{r=0}^{j} (-1)^{j-r} \cdot \sum_{s=0}^{i} k^{i-s} \cdot \binom{i}{s} \cdot  
\sum_{l=0}^{r} (-1)^{r-l} \binom{r}{l} \cdot l^s=\\
&\sum_{r=0}^{j} \sum_{l=0}^{r} (-1)^{j-l} \cdot \binom{r}{l} \cdot (k+l)^{i}=\\
&\sum_{l=1}^{j+1} (-1)^{j+1-l} \cdot (k-1+l)^{i} \cdot \binom{j+1}{l}=
b_{i,j+1,k-1}+(-1)^j \cdot (k-1)^i.  \ \square{}
\end{align*}

The following corollary shows four interesting special cases.
\begin{cor} \ 
Under the assumptions of the preceding theorem, the following identities are valid:
\begin{align}
\text{a) }  &\sum_{r=0}^{j} (-1)^{j-r} \cdot b_{i,r,1} = b_{i,j+1,0}, 
				\text{ if } i \ge 1 \text{ and } j \ge 1, \label{a25}\\ 
\text{b) }  &\sum_{r=0}^{i-1} (-1)^{i-1-r} \cdot b_{i,r,k} = i! + (-1)^{i-1} \cdot (k-1)^i,  		\label{a26}\\
\text{c) }  &\sum_{r=0}^{i} (-1)^{r} \cdot b_{i,r,k} = (k-1)^i,  		\label{a27}\\
\text{d) }  &\sum_{r=0}^{i} (-1)^{i-r} \cdot b_{i,r,0} =1. \label{a28}
\end{align}
\end{cor}  
\textit{Proof.} a), b) and c) follow from (\ref{a24}) for $k=1$, $j=i$ resp. $j=i$. d) follows for $k=0$ from c) 	$ \ \square{} $

We continue with summing over the second parameter.
\begin {theorem} 
For integers $i \ge 0$, $i \ge j \ge 0$ and $k \ge 1-j$
\begin{align}
&\text{a) } \sum_{r=j}^{i} (-1)^{i-r} \cdot \binom{r+k-1}{r-j} \cdot b_{i,r,k}=b_{i,j,0}.\label{a29}\\ 
&\text{b) } \sum_{r=j}^{i} (-1)^{i-r} \cdot \binom{r+k-1}{r-j} \cdot b_{i,r,0}=b_{i,j,k}.\label{a29a} 
\end{align}
\end{theorem}
\textit{Proof.}  The proof of a) is performed by induction on $i$. For $i=j$ both sides have according to (\ref{a14}) the value $j!$. If (\ref{a29}) holds for $i$ as induction hypothesis, then due to (\ref{a9}) and well known identities about the binomial coefficients (cf. for example \cite{DLMF})
\begin{align*}
	&\sum_{r=j}^{i+1} (-1)^{i+1-r} \cdot \binom{r+k-1}{r-j} \cdot b_{i+1,r,k}=\\
	&\sum_{r=j}^{i+1} (-1)^{i+1-r} \cdot \binom{r+k-1}{r-j} \cdot r \cdot b_{i,r-1,k}+\\
	&\qquad \qquad \qquad \sum_{r=j}^{i} (-1)^{i+1-r} \cdot \binom{r+k-1}{r-j} 
	\cdot (r+k) \cdot b_{i,r,k}=\\
	&j \cdot \sum_{r=j-1}^{i} (-1)^{i-r} \cdot b_{i,r,k} \cdot \binom{r+k}{r-j+1}=\\
	&j \cdot \sum_{r=j-1}^{i} (-1)^{i-r} \cdot b_{i,r,k} \cdot \binom{r+k-1}{r-j}+
	j \sum_{r=j-1}^{i} (-1)^{i-r} b_{i,r,k} \binom{r+k-1}{r-(j-1)}=\\
	&j \cdot (b_{i,j,0}+b_{i,j-1,0})=b_{i+1,j,0}.
\end{align*}
b) follows from a) by means of the inversion formula for binomial coefficients.
$\square{}$

For $k=1-j$ and $k=2-j$ interesting special cases result from the preceding theorem, which are omitted here for brevity. 

The following theorem is useful for the computation of \textit{MSNs} with high values of the third parameter.
\begin {theorem} \label{Interpretation}
For integers $i \ge 0$, $j \ge 0$, $k_2 \ge 0$, $k \ge 0$, $l > 1$ and $k_1$
\begin{align}
&\text{a) }
b_{i,j,k_1+k_2}=\sum_{r=0}^{k_2} \binom{k_2}{r} \cdot b_{i,j+r,k_1},\label{a30} \\
&\text{b) }
b_{i,j,l \cdot k}=\sum_{r=0}^{min(i-j,(l-1) \cdot k)} \binom{(l-1) \cdot k}{r} \cdot b_{i,j+r,k}.\label{n30}
\end{align}
\end{theorem}
\textit{Proof.} a) Due to (\ref{a1}), well known identities about the binomial coefficients (cf. for example \cite{DLMF}) and index transformation $r=s-k_2$
\begin{align*}
	&\sum_{r=0}^{k_2} \binom{k_2}{r} \cdot b_{i,j+r,k_1}=
	\sum_{s=0}^{j+k_2} (s+k_1)^i \cdot (-1)^{j-s} \sum_{r=0}^{k_2} \binom{k_2}{r} 
	\cdot \binom{j+r}{s} \cdot (-1)^{r} =\\
	&\sum_{r=0}^{j} (r+k_1+k_2)^i \cdot (-1)^{j-r} \cdot \binom{j}{r}=b_{i,j,k_1+k_2}. 	
\end{align*}
b) follows from a) due to $l \cdot k =(l-1) \cdot k +k$ and (\ref{a13}). $\square{}$

For $k_1=0$, the preceding theorem shows due to (\ref{a5}) a tight connection between $b_{i,j,k_2}$ and a sum of "`ordinary"' Stirling numbers of the second kind. 

Further, for $k_1=0$ and $k_2=k \ge 0$, (\ref{a30}) allows a combinatorial interpretation of $b_{i,j,k}$: If $M_{j+r}$, $r=0, \ldots, k$, are subsets of $M_{j+k}$ with $|M_{j+r}|=j+r$, $r=0, \ldots,k$, and analogously for another set $M_{i}$ with $|M_{i}|=i$, then $b_{i,j,k}$ is the number of functions $f: M_i \rightarrow M_{j+k}$ with $M_j \subset f(M_i)$, because $b_{i,j+r,0}=S_{i,j+r} \cdot (j+r)!$ is the number of surjections from $M_{i}$ onto $M_{j+r}$ and there are $\binom{k}{r}$ subsets $M$ with $M_j \subset M \subset M_{j+k}$ and $|M|=j+r$.  

After the summation over the first and the second parameter of the \textit{MSNs}, we now turn to the summation over the third parameter. The following theorem is an inversion of (\ref{a30}).
\begin {theorem} 
For integers $i \ge 0$, $j \ge 0$, $k_2 \ge 0$ and real $k_1$
\begin{equation}
b_{i,j+k_2,k_1}=\sum_{r=0}^{k_2} \binom{k_2}{r} \cdot (-1)^{k_2-r} \cdot b_{i,j,k_1+r}.\label{a31} 
\end{equation}
\end{theorem}
\textit{Proof.} Due to (\ref{a30})  
	\[b_{i,j,k_1+r}=\sum_{s=0}^{r} \binom{r}{s} \cdot b_{i,j+s,k_1}
	\]
and the inversion formula for binomial coefficients, we get
\begin{align*}
	&\sum_{r=0}^{k_2} (-1)^{k_2-r} \cdot \binom{k_2}{r} \cdot b_{i,j,k_1+r}=
	\sum_{r=0}^{k_2} (-1)^{k_2-r} \cdot \binom{k_2}{r} \cdot \sum_{s=0}^{r} \binom{r}{s} \cdot b_{i,j+s,k_1}=\\
	&\sum_{s=0}^{k_2} (-1)^{k_2-s} \sum_{r=s}^{k_2} (-1)^{s-r} \cdot \binom{k_2}{r} \cdot \binom{r}{s} \cdot b_{i,j+s,k_1}= b_{i,j+k_2,k_1}. \ \square{}
\end{align*}

If $k_1=0$, then for $k_2=i-j$ and $k_2=i$ the preceding theorem contains due to (\ref{a14}) and (\ref{a13}) two interesting special cases:
\begin{cor} For integers $j \ge 0$ and $i \ge j$   
\begin{align}
& a) &\sum_{r=0}^{i-j} \binom{i-j}{r} \cdot (-1)^{i-j-r} \cdot b_{i,j,r}=i! , \ \label{a32} \\
& b) &\sum_{r=0}^{i} \binom{i}{r} \cdot (-1)^{i-r} \cdot b_{i,j,r}=0, \ 
\text{ if } j>0. \ \label{a33}
\end{align} 
\end{cor}  
\begin {theorem} 
For integers $i \ge 0$, $j \ge 0$, and $k \ge 0$
\begin{equation}
b_{i,j,k+1}=b_{i,j,1}+\sum_{r=1}^{k} b_{i,j+1,r}.\label{a34} 
\end{equation}
\end{theorem}
\textit{Proof.} 
Induction on $k$ by means of (\ref{a8}). $\square{}$

\begin{cor} For integers $i \ge 1$, $l \ge 1$ and $k \ge 1$   
\begin{align}
& a) &\sum_{r=0}^{k-1} b_{i,1,r}&= k^i, \label{k_i} \\
& b) &\sum_{r=0}^{k-1} b_{i,1,l+r}&= (k+l)^i -l^i. \label{k_i_l}
\end{align} 
\end{cor} 
\textit{Proof.} a) follows from (\ref{a34}) for $j=0$, because $b_{i,0,k}=k^i$ and $b_{i,0,1}=0$ (cf. (\ref{a4})). b) follows from a), if the sum range $\{0, \ldots, k+l-1\}$ is split into the disjoint sets $\{0, \ldots, l-1\}$ and $\{l, \ldots, k+l-1\}$. \ \ $\square{}$

\paragraph{}
For completeness, also the inverses of the \textit{MSNs} are investigated. Because the Stirling numbers $s_{i,j}, i\ge 0, j \ge 0$, of the first kind are inverses of the Stirling numbers  of the second kind, they also play here an important role. Since there are different definitions, we repeat their (recursive) definition in order to guard against misunderstandings:
\begin{equation}
s_{i+1,j}=s_{i,j-1}-i \cdot s_{i,j} \label{sn1}
\end{equation}
with $s_{0,j}=s_{i,0}=0$,  for $i>0, j>0$ and $s_{0,0}=1$ (cf. e.g. \cite{Heise}, \cite{DLMF} or \cite{Abramowitz}). We generalize them to the \textbf{\textit{moment generating Stirling numbers of the first kind, "`MSN1"'}} for integers $i\ge 0,\ j\ge 0$ and real $k$ by setting
\begin{equation}
	c_{i,j,k}=\sum_{r=j}^{i} \binom{r}{j} \cdot (-k)^{r-j} \cdot s_{i,r}.
	\label{sn2}
\end{equation}
Because the \textit{MSN1s} are not in the focus of this paper, we do not list their properties and connections between them and the \textit{MSNs} in detail, but only remark that $c_{i,j,0}=s_{i,j}$.
\begin {theorem} 
For integers $i \ge 0$, $j \ge 0$ and real $k_1$ and $k_2$
\begin{equation}
\sum_{r =j}^{i} \frac{1}{r!} \cdot b_{i,r,k_1} \cdot c_{r,j,k_2}= 
	\binom{i}{j} \cdot (k_1 - k_2)^{i-j}. \label{a46} 
\end{equation}
\end{theorem}
\textit{Proof.} First, we show (\ref{a46}) for $k_1=k_2=k$: Due to (cf. e.g. \cite{Heise}) 
	\[\sum_{r=j}^{i} S_{i,r} \cdot s_{r,j}= \delta_{i,j}, \]
(\ref{a20}) and (\ref{a5}), the assertion follows from 
\begin{align*}
	&\sum_{r=j} ^{i} \frac{b_{i,r,k}}{r!} \cdot c_{r,j,k}= 
	\sum_{r=j} ^{i} \frac{b_{i,r,k}}{r!} \cdot \sum_{l = j}^{r} \binom{l}{j} \cdot (-k)^{l-j} \cdot s_{r,l} =\\
	&=\sum_{l=j}^{i} \binom{l}{j} \cdot (-k)^{l-j} \cdot \sum_{r=j}^{i} \frac{1}{r!} \cdot 
	\sum_{s=r}^{i} \binom{i}{s} \cdot b_{s,r,0} \cdot k^{i-s} \cdot s_{r,l} =\\
	&=\sum_{l=0}^{i} \binom{l}{j} \cdot (-k)^{l-j} \cdot \binom{i}{l} \cdot k^{i-l} = 
	\binom{i}{j} \cdot k^{i-j} \cdot \sum_{l=0}^{i-j} \binom{i-j}{l} \cdot (-1)^l=
	\delta_{i,j} 
\end{align*}
and the general case from
\begin{align*}
	&\sum_{r=j} ^{i} \frac{b_{i,r,k_1}}{r!} \cdot c_{r,j,k_2}= 
	\sum_{r=j} ^{i} \frac{1}{r!} \cdot \sum_{s=0} ^{i} \binom{i}{s} \cdot (k_1-k_2)^{i-s} \cdot
	b_{s,r,k_2} \cdot c_{r,j,k_2} = \\
	&=\binom{i}{j} \cdot (k_1-k_2)^{i-j}. 	\ \square{}
\end{align*}
\subsection{Advanced Topics}
After we've listed the basics, let's look at some advanced topics which can be derived from the results of the preceding subsection. First, we represent the MSNs
(of the second kind) by the multinomial coefficients (cf. \cite{Abramowitz}).

\begin {theorem} 
For integers $i \ge 0$, $l \ge 1$, $j_1 \ge 0,\ldots,j_l \ge 0$ and $k_1 ,\  \ldots,k_l $
\begin{equation}
b_{i,j_1+\ldots+j_l,k_1+\ldots+k_l}=\sum_{} \binom{i}{i_1, \ldots , i_l} \cdot 
\prod _{r=1}^{l} b_{i_r,j_r,k_r}, \label{a36} \\
\end{equation}
where the sum is over all $i_1 \ge 0, \ldots , i_l \ge 0$ with $i=i_1 +\ldots +i_l$.
(Of course, the sum can be restricted to the values $i_r \ge j_r$, $r=1,\ldots, l$, due to (\ref{a13})). 
\end{theorem}
\textit{Proof.} Induction on $l$: For $l=1$, the assertion is obvious. If it holds for $l$, then due to (\ref{a17})
\begin{align*}
	&\sum_{} \binom{i}{i_1, \ldots , i_{l+1}} \cdot \prod _{r=1} ^{l+1} b_{i_r,j_r,k_r} =\\
	&\sum_{i_{l+1}=0} ^{i} \binom{i}{i_{l+1}} \cdot b_{i_{l+1},j_{l+1},k_{l+1}} \cdot 
	\sum_{} \binom{i-i_{l+1}}{i_1, \ldots , i_{l}} \cdot \prod _{r=1} ^{l} b_{i_r,j_r,k_r}=\\
	&\sum_{i_{l+1}=0} ^{i} \binom{i}{i_{l+1}} \cdot b_{i_{l+1},j_{l+1},k_{l+1}} \cdot b_{i-i_{l+1},j_1+\ldots +j_l,k_1+ \ldots + k_l}=\\
	&b_{i,j_1+\ldots+j_{l+1},k_1+\ldots+k_{l+1}}. \ \square{}
\end{align*}

Obviously, (\ref{a17}) is a special case for $l=2$.
\begin{cor}  
For integers $i \ge 0$, $l \ge 1$, $j \ge 0,j_1 \ge 0,\ldots,j_l \ge 0$  and $k \ge 0$
\begin{align}
\text{a) }  &b_{i,j_1+\ldots+j_l,k}=\sum_{} \binom{i}{i_1, \ldots , i_{l+k}} \cdot 
\prod _{r=1}^{l} b_{i_r,j_r,0}, \label{a37} \\
\text{b) }  &b_{i,j,k}=\sum_{i_1>0, \ldots, i_j>0} \binom{i}{i_1, \ldots , i_{j+k}}. \label{a38} 
\end{align}
\end{cor}  
\textit{Proof.} a) Due to (\ref{a36}) with $s=i_{l+1}+\ldots+i_{l+k}=i-(i_1+\ldots+i_l)$ and (\ref{a20}),  
\begin{align*}
	&\sum_{} \binom{i}{i_1, \ldots , i_{l+k}} \cdot \prod _{r=1}^{l} b_{i_r,j_r,0}=\\
	&\sum_{s=0}^{i} \binom{i}{s} \cdot \sum_{} \binom{s}{i_{l+1}, \ldots , i_{l+k}} \cdot 
	\sum_{} \binom{i-s}{i_{1}, \ldots , i_{l}} \cdot \prod _{r=1}^{l} b_{i_r,j_r,0}=b_{i,j_1+\ldots +j_l,k}.
\end{align*}

b) follows from (\ref{a37}) for $l=j$, $j_1=\ldots=j_l=1$, because $b_{i_r, 1,0}=0$, if $i_r=0$, $r=1, \ldots, j$ and $b_{i_r, 1,0}=1$, if $i_r>0$ (cf. (\ref{a6})). \ $\square{}$

One should note that the formulae 
	\[j^i=\sum \binom{i}{i_1, \ldots, i_j} \]
and (\ref{a38}) for $k=0$ 
	\[b_{i,j,0}=\sum_{i_1 >0, \ldots, i_j>0} \binom{i}{i_1, \ldots, i_j}	\]
are very similar, but not equal. Additionally, (\ref{a38}) is due to (\ref{a5}) for $k=0$ a representation of $S_{i,j}\cdot j!$, and the proof of (\ref{a37}) shows that it generalizes (\ref{a20}) for $k \ge 0$.

As last topic in this section, we deal with some generating functions, which allow to develop and proof further theorems about the MSNs. They generalize some results from \cite{Abramowitz}. 
\begin {theorem} 
For integers $i \ge 0$, $j \ge 0$, and real $k$		

a) The (ordinary) generating function
\begin{equation}
F_{j,k}(x)= \sum_{i \ge 0} b_{i,j,k} \cdot x^i = \frac{1}{1-k \cdot x} \cdot \frac{j! \cdot x^j}{\prod_{r=1}^{j}\left(1-(k+r)\cdot x \right)} , \label{a39} 
\end{equation}
b) the exponential generating function for integer $k$
\begin{equation}
E_{j,k}(x)=\sum_{i \ge 0} b_{i,j,k} \cdot \frac{x^i}{i!} = \left(e^x-1\right)^j \cdot e^{k \cdot x}, \label{a40} \\
\end{equation}
c) the double exponential generating function for integer $k$
\begin{equation}
E_{j}(x,y)=\sum_{i \ge 0, k \ge 0} b_{i,j,k} \cdot \frac{x^i}{i!} \cdot \frac{y^k}{k!} =
\left(e^x-1\right)^j \cdot exp(e^x \cdot y), \label{a41} \\
\end{equation}
d) for integer $k$
\begin{equation}
E_{k}(x,z)=\sum_{i \ge 0, j \ge 0} b_{i,j,k} \cdot \frac{x^i}{i!} \cdot \frac{z^j}{j!} =
e^{k \cdot z}  \cdot exp((e^x -1) \cdot z), \label{a42} \\
\end{equation}
e) the binomial generating function
\begin{equation}
B_{i,k}(x)=\sum_{j \ge 0} b_{i,j,k} \cdot \binom{x}{j}= x^i \cdot (k+1)^i, \label{a43} \\
\end{equation}
f) the mixed exponential binomial generating function 
\begin{equation}
EB_{k}(x,y)=\sum_{i \ge 0, j \ge 0} b_{i,j,k} \cdot \frac{y^i}{i!} \cdot \binom{x}{j}= e^{(k+1)\cdot x \cdot y}. \label{a44} \\
\end{equation}
\end{theorem} 
\textit{Proof.} a) For $j=0$, $F_{0,k}(x)=\frac{1}{1-k \cdot x}$. Due to (\ref{a9}), 
\begin{align*}
	&F_{j+1,k}(x)= \sum_{i \ge 0} b_{i,j+1,k} \cdot x^i =
	(j+1) \cdot x \cdot F_{j,k}(x) +(j+1+k) \cdot x \cdot F_{j+1,k}(x)
\end{align*}
and therefore
	\[F_{j+1,k}(x)=\frac{(j+1) \cdot x }{1-(j+1+k) \cdot x} \cdot F_{j,k}(x).\]
The assertion follows by induction on $j$. b) The exponential generating function of the Stirling numbers of the second kind is (cf. \cite{Abramowitz})
	\[\sum_{i \ge 0} S_{i,j} \cdot \frac{x^i}{i!} = \frac{(e^x-1)^j}{j!}.\]
Therefore, due to (\ref{a30})	
\begin{align*}
	&E_{j,k}(x)=\sum_{i \ge 0} b_{i,j,k} \cdot \frac{x^i}{i!} = 
	\sum_{i \ge 0}\sum_{r = 0}^{k} \binom{k}{r} \cdot b_{i,j+r,0} \cdot \frac{x^i}{i!}=
		\left(e^x-1\right)^j \cdot e^{k\cdot x}.
\end{align*}
c) and d) follow from (\ref{a40}).
e) Due to (\ref{a20}) and (\ref{a7})
\begin{align*}
B_{i,k}(x)=&\sum_{j = 0}^{i} b_{i,j,k} \cdot \binom{x}{j}= 
\sum_{j = 0} ^{i} \binom{x}{j} \sum_{r = 0}^{i} \binom{i}{r} \cdot k^{i-r} \cdot b_{r,j,0}= x^i \cdot (k+1)^i.
\end{align*}
f) follows from e). $\square{}$

After it has been shown so far, how the \textit{MSNs} generalize the Stirling numbers of the second kind in a natural way and that they have more properties than are known for the generalized Stirling numbers in \cite{Hsu}, applications will be presented in the next section.

\section{Applications} \ 
{\allowdisplaybreaks

The purpose of this section is to show that the moments of the discrete phase type distribution ($PH_{D}$), the $PH_{D}$-renewal times, the recurrence time of a discrete time Markov chain and the generalized and ordinary negative binomial distribution (with its special case geometric distribution) can be expressed by the \textit{MSNs} as closed formulae. This is possible by means of the framework developed in \cite{Fra} and \cite{Fra2}, which generalizes all these distributions: 

Let (\(X_r,\ r \geq 0\)) be a Discrete Time Markov Chain (DTMC) on
the space \(S = \{1, 2, \ldots\, m+1\}\) with transition probability matrix
\(\textbf{\textit{P}}\) (cf. \cite{Kulk}) and subset \(M \subset S,\
M \neq \emptyset, \ \overline{M} \neq \emptyset\). As in \cite{Fra} we subdivide \(\textbf{\textit{P}}=(p_{ij})\), \(i, j \in S\), in the following way into submatrices:
{\allowdisplaybreaks
\begin{equation}\label{1}
    \textbf{\textit{P}} =
\begin{pmatrix}
    \textbf{\textit{P}}_{M} & \textbf{\textit{P}}_{M\overline{M}} \\
    \textbf{\textit{P}}_{\overline{M}M} & \textbf{\textit{P}}_{\overline{M}} \\
\end{pmatrix}
\end{equation}
where \(\textbf{\textit{P}}_{M} = (p_{ij})\) with \(i, j \in
M, \textbf{\textit{P}}_{M\overline{M}} = (p_{ij})\) with \(i \in
M, \ j \in \overline{M}\),
\(\textbf{\textit{P}}_{\overline{M}M} = (p_{ij})\) with
\(i \in \overline{M},\ j \in M\) and \(\textbf{\textit{P}}_{\overline{M}} = (p_{ij})\)
with $i, j \in \overline{M}$, and assume that $\textbf{\textit{I}}-\textbf{\textit{P}}_{M}$ is invertible.

Accordingly, we also subdivide the \(|S|\)-dimensional column vector \(\textbf{\textit{e}}\) that consists of exclusively 1 in all components of \(S\), into the column vectors \(\textbf{\textit{e}}_M\) and \(\textbf{\textit{e}}_{\overline{M}}\).

For any integer \(k \geq 1\), we define the \(|M|\times |\overline{M}|\)-matrix random variable $\mathcal{N}_{k}=(\mathcal{N}_{k,i,j})_{i \in M, j \in \overline{M}}$, which describes the time $r$ until the number of visits in \(\overline{M}\) is \(k\) for the first time and $X_r=j \in \overline{M}$, starting in \(i \in M\). 
\(\mathcal{R}_{k}=(\mathcal{R}_{k,i,j})_{i \in M,j \in M}\) is the \(|M|\times |M|\)-matrix random variable, which describes the time $r$ until the number of returns to \(M\) is \(k\) for the first time and $X_r=j \in M$, starting in \(i \in M\). 
\(\mathcal{N}_{k}\) is the \textit{\(k\)-th passage time into} \(\overline{M}\) and 
\(\mathcal{R}_{k}\) the \textit{\(k\)-th recurrence time of \(M\)}.  

If \(|M|=1\), \(\mathcal{R}_{1}\) is the ordinary recurrence time (cf. \cite{Kulk}). If in addition \(|\overline{M}|=1\) and 
\begin{equation}
    \textbf{\textit{P}} =
\begin{pmatrix}
    1-p & p \\
    1-q & q \\
\end{pmatrix},
\end{equation}
$0 \le p \le 1, 0 \le q \le 1$, \(\mathcal{N}_{1}\) and \(\mathcal{R}_{1}\) are in \cite{Gosh} and \cite{Triv} used to model software reliability. They describe the time to the next failure resp. the time between failures. In this case the distribution of \(\mathcal{N}_{k}\) is called \textit{alternating negative binomial distribution}, because it counts in a series of independent Bernoulli trials the number of trials until the $k-$th success occurs, where the success probability alternates between $p$ and $q$ depending on whether the preceding trial was a failure (or the first one) or a success. If in addition $p=q$ the distribution of \(\mathcal{N}_{k}\) is the ordinary negative binomial distribution. Because of this special case, we call the distribution of \(\mathcal{N}_{k}\) \textit{generalized negative binomial distribution} and for \(k=1\) the \textit{generalized geometric distribution}.

We can also interchange the roles of \(M\) and \(\overline{M}\) and start the Markov chain in \(i \in \overline{M}\). Then, we define in an analogous way for integer \(k>0\) the matrix random variable 
$\overline{\mathcal{N}}_{k}$, which counts the visits in $M$, and $\mathcal{\overline{R}}_{k}$, which counts the returns to $\overline{M}$.

In \cite{Fra} and \cite{Fra2}, recursive formulae for all $P(\mathcal{R}_{k}=n)$ and $P(\mathcal{N}_{k}=n)$, \(n \geq 0\) and \(k \geq 0\), are deduced, but the cases \(P(\mathcal{R}_{1}=n)\) and \(P(\mathcal{N}_{1}=n)\) are of special interest, because they allow the simplest formulae and are the base for the computation of all values for \(k > 1\). 
The following Lemma allows to compute the (conditional) probabilities of the first return to \(M\) and (first) passage time into \(\overline{M}\) for all starting states \(i \in M\) and all end states \(j \in M\) resp. \(j \in \overline{M}\). It also shows the connection between the recurrence and the first passage time and establishes the relation to discrete phase type distributions and phase type renewal processes. 
\begin{lemma} \label{lemma1_neu}
For \(n \geq 1\) the following formulae hold:
\begin{align}
  \text{a) }  &P(\mathcal{N}_{1}=n) & = & \textbf{\textit{P}}^{n-1}_{M} \cdot \textbf{\textit{P}}_{M\overline{M}}, \label{lemma1a}\\
  \text{b) } &P(\mathcal{R}_{1}=n) & = & \textbf{\textit{P}}_{M\overline{M}} \cdot \textbf{\textit{P}}^{n-2}_{\overline{M}} \cdot \textbf{\textit{P}}_{\overline{M}M} \text{ for } n \geq 2,  \label{lemma1da} \nonumber \\
	   \ \ &P(\mathcal{R}_{1}=1) & = & \textbf{\textit{P}}_{M},\\ 	
	\text{c) } &P(\mathcal{R}_{1}=n) & = & \textbf{\textit{P}}_{M\overline{M}} \cdot P(\overline{\mathcal{N}}_{1}=n-1), \label{lemma1c}\\
	  \text{d) } &P(\mathcal{N}_{1}=n) & = & \textbf{\textit{P}}_{M}^{n-1} \cdot (\textbf{\textit{I}}-\textbf{\textit{P}}_{M})\cdot \textbf{\textit{e}}_M \text{ , if } |\overline{M}|=1.
		\label{lemma1d}
\end{align}
\end{lemma}

A detailed proof can be found in \cite{Fra2} or for a) and d) in \cite{Neuts} or \cite{Neuts_geom}.
Analogous formulae hold for $P(\overline{\mathcal{N}}_{1}=n)$ and $P(\overline{\mathcal{R}}_{1}=n)$.

The preceding lemma shows that for $|\overline{M}|=1$ $\overline{\mathcal{R}}_{1}$ is 
a shifted $PH_D-$ distribution, i.e., if a random variable $X$ on $\{1,2,\ldots,\}$ is defined by $P(X=n)=P(\overline{\mathcal{R}}_1=n+1)$, $n \ge 1$, then $X$ is $PH_D(\textbf{\textit{P}}_{\overline{M}M},\textbf{\textit{P}}_{M})$- distributed with matrix $\textbf{\textit{P}}_{M}$ and initial vector $\textbf{\textit{P}}_{\overline{M}M}$, because 
$P(X=n)=P(\overline{\mathcal{R}}_{1}=n+1)=\textbf{\textit{P}}_{\overline{M}M} \cdot \textbf{\textit{P}}^{n-1}_{M} \cdot (\textbf{\textit{I}}- \textbf{\textit{P}}_{M}) \cdot \textbf{\textit{e}}_{M}$. 

If $\overline{M}=\{m+1\}$ and $\textbf{\textit{P}}_{\overline{M}}=(0)$, we get a discrete time phase type renewal process, because after entering state $m+1$ the process starts immediately again with the same initial condition $\textbf{\textit{P}}_{\overline{M}M}$ (cf. \cite{Kulk}). The visits in state $m+1$ are the renewal points and the number of renewals at time $n \ge 0$ is a renewal process. In this case $\overline{\mathcal{R}}_k$ is the time until the $k$-th renewal. 
The moments (and especially the expected number of renewals at time $n \in \mathbb{N}$) of the renewal process were computed in \cite{Fra}.  

Formulas as in (\ref{lemma1a}) also occur in discrete time Markovian arrival processes (DMAP), cf. e.g. \cite{Briem} or \cite{Tran}, if $|M|=|\overline{M}|$. So, some of our results can also be applied to DMAPs. 

In the following definitions of the moments, we assume, as usual in this context, \(0^{0}=1\). 
For integers \(k \geq 1\) and \(m \geq 0\)
\begin{align} 
M_{m}(\mathcal{N}_{k})= &\sum_{n=0}^{\infty} n^{m} \cdot P(\mathcal{N}_{k}=n), \label{MmNk}\\
M_{m}(\mathcal{R}_{k})= &\sum_{n=0}^{\infty} n^{m} \cdot P(\mathcal{R}_{k}=n) \label{MmRk}
\end{align}
and analogously $M_{m}(\overline{\mathcal{N}}_{k})$ and  $M_{m}(\overline{\mathcal{R}}_{k})$.
In this section we focus on the computation of $M_{m}(\mathcal{N}_{k})$, $M_{m}(\mathcal{R}_{k})$ and $M_{m}(\overline{\mathcal{R}}_{k})$ 
 
The following lemma shows recursive formulae for the moments of $\mathcal{N}_{1}$ and $\mathcal{R}_{1}$, which are the basis for the closed formulae in the sequel.

\begin{lemma} \label{lemma2_neu}
Under the assumptions of this section
\begin{align}
  \text{a) }  &M_0(\mathcal{N}_{1}) & = & (\textbf{\textit{I}}-\textbf{\textit{P}}_M)^{-1} \cdot \textbf{\textit{P}}_{M\overline{M}}, \label{lemma2a}\\
  \text{b) } &M_m(\mathcal{N}_{1}) & = & (\textbf{\textit{I}}-\textbf{\textit{P}}_M)^{-1} \cdot \left(\textbf{\textit{P}}_{M\overline{M}} + \textbf{\textit{P}}_M \cdot \sum_{j=0}^{m-1} \binom{m}{j} \cdot  M_j(\mathcal{N}_{1})\right)   \label{lemma2b} \\
	& & & \qquad \qquad \text{ for } m \geq 1, \nonumber \\
	\text{c) } &M_m(\mathcal{R}_{1}) & = & \textbf{\textit{P}}_{M} + 
	\textbf{\textit{P}}_{M\overline{M}} \cdot \sum_{j=0}^{m} \binom{m}{j} \cdot  
	M_j(\overline{\mathcal{N}}_{1})   \text{ for } m \geq 0. \label{lemma2c}       
\end{align}
\end{lemma}
\textit{Proof.} a) follows from (\ref{lemma1a}), b) is a consequence of $n^m=\sum_{j=0}^{m}\binom{m}{j} \cdot (n-1)^j$ and c) is based on (\ref{lemma1c}).  $\ \square{} $

Detailed proofs can be found under more general conditions in \cite{Fra2}, Theorem 3.3 and 3.5. The following theorem shows the first benefits of the \textit{MSNs}.

\begin{theorem} For $m \geq 0$
\begin{align}
  \text{a) }  &  M_m(\mathcal{N}_{1})  =  
	\sum_{j=0}^{m} b_{m,j,1} \cdot \textbf{\textit{P}}_M^{j} \cdot 
	(\textbf{\textit{I}}-\textbf{\textit{P}}_M)^{-j-1} \cdot \textbf{\textit{P}}_{M\overline{M}}, \label{n3_14} \\
	\text{b) }	&M_m(\mathcal{R}_{1})  =   
	\textbf{\textit{P}}_{M} + \textbf{\textit{P}}_{M\overline{M}} \cdot  
	\sum_{j=0}^{m} b_{m,j,2} \cdot \textbf{\textit{P}}_{\overline{M}}^{j} \cdot (\textbf{\textit{I}}-\textbf{\textit{P}}_{\overline{M}})^{-j-1} \cdot \textbf{\textit{P}}_{\overline{M}M}.
	\label{n3_15}
\end{align}
\end{theorem}
\textit{Proof.} a) Induction on $m$: For $m=0$ the assertion follows from (\ref{lemma2a}) due to $b_{0,0,1}=1$. If (\ref{n3_14}) holds for $m$ as induction hypothesis, then due to (\ref{lemma2b}), (\ref{a10}) and $b_{m+1,0,1}=1$
\begin{align*}
	&(\textbf{\textit{I}}-\textbf{\textit{P}}_M) \cdot M_{m+1}(\mathcal{N}_{1})  =
	\textbf{\textit{P}}_{M\overline{M}} + \textbf{\textit{P}}_M \cdot \sum_{j=0}^{m} \binom{m+1}{j} \cdot  M_j(\mathcal{N}_{1})= \\
	&=\textbf{\textit{P}}_{M\overline{M}} + \textbf{\textit{P}}_M \cdot 
	\sum_{j=0}^{m} \binom{m+1}{j} \cdot
	 \sum_{r=0}^{j} b_{j,r,1} \cdot  
	\textbf{\textit{P}}_M^{r} \cdot (\textbf{\textit{I}}-\textbf{\textit{P}}_M)^{-r-1} \cdot \textbf{\textit{P}}_{M\overline{M}} =\\
	&=\textbf{\textit{P}}_{M\overline{M}} + \sum_{r=0}^{m} b_{m+1,r+1,1} \cdot \textbf{\textit{P}}_M^{r+1} (\textbf{\textit{I}}-\textbf{\textit{P}}_M)^{-r-1} \cdot \textbf{\textit{P}}_{M\overline{M}}= \\
	&=\sum_{r=0}^{m+1} b_{m+1,r,1} \cdot \textbf{\textit{P}}_M^{r} (\textbf{\textit{I}}-\textbf{\textit{P}}_M)^{-r} \cdot \textbf{\textit{P}}_{M\overline{M}}. 
\end{align*}
b) follows from a) due to (\ref{lemma2c}) and (\ref{a11}). 	\	$\square{}$

An alternative proof of (\ref{n3_14}) is possible by means of the factorial moments (cf. \cite{Neuts_geom}) $F_m(\mathcal{N}_1)=m ! \cdot \textbf{\textit{P}}_M^{m-1} \cdot (\textbf{\textit{I}}-\textbf{\textit{P}}_M)^{-m-1} \cdot \textbf{\textit{P}}_{M\overline{M}}$, $M_m(\mathcal{N}_1)=\sum_{j=0}^{m} \frac{1}{j!} \cdot b_{m,j,0} \cdot F_j(\mathcal{N}_1)$ and a tricky computation.

The following lemma allows the iterative computation of $M_m(\mathcal{N}_k)$ and $M_m(\mathcal{R}_k)$ without additional assumptions and is the basis of the closed formulas in the sequel.

\begin{lemma}  For $m \ge 0$ and $k \ge 2$ 
\begin{align}
  \text{a) }  & M_m(\mathcal{R}_k) = 
	     \sum_{j=0}^{m} \binom{m}{j }\cdot M_j(\mathcal{R}_{k-1}) \cdot M_{m-j}(\mathcal{R}_{1}),     
			\label{n3_16} \\
	\text{b) }  & M_m(\mathcal{N}_k) = 
	     \sum_{j=0}^{m} \binom{m}{j }\cdot M_{m-j}(\mathcal{N}_{1}) \cdot M_{j}(\mathcal{\overline{R}}_{k-1}).  \label{n3_17}
\end{align}	
\end{lemma}

The proof follows from $\mathcal{R}_{k}=\mathcal{R}_{k-1}+\mathcal{R}_{1}$ and $\mathcal{N}_{k}=\mathcal{N}_{1}+\overline{\mathcal{R}}_{k-1}$. Detailed proofs are performed in \cite{Fra2}, Theorems 3.8 and 3.10.

We call the stochastic matrix $\textbf{\textit{P}}$ \textbf{\textit{M-commutable}}, if for all integers $r \ge 0$ and $s \ge 0$ 
$\textbf{\textit{P}}_{M\overline{M}} \cdot \textbf{\textit{P}}_{\overline{M}}^s \cdot \textbf{\textit{P}}_{\overline{M} M} \cdot \textbf{\textit{P}}_{M}^r=\textbf{\textit{P}}_{M}^r \cdot \textbf{\textit{P}}_{M\overline{M}} \cdot \textbf{\textit{P}}_{\overline{M}}^s \cdot \textbf{\textit{P}}_{\overline{M}M}$. 
Analogously, \textbf{\textit{$\overline{M}$-commutability}} is defined. 

For example, if $|M|=1$, $\textbf{\textit{P}}_{M}= p \cdot \textbf{\textit{I}}$ with $p \ge 0$, $\textbf{\textit{P}}_{M\overline{M}}^T=p \cdot \textbf{\textit{P}}_{\overline{M}M}$ or $\textbf{\textit{P}}$ consists of just the same lines, then $\textbf{\textit{P}}$ is $M-$commutable.

Under the additional assumption of commutability or if all lines  of $\textbf{\textit{P}}_{\overline{M}}$ (resp. $\textbf{\textit{P}}_{M}$) have equal sums (i.e. 
$\textbf{\textit{P}}_{\overline{M}} \cdot \textbf{\textit{e}}_{\overline{M}}=s_{\overline{M}} \cdot \textbf{\textit{e}}_{\overline{M}}$), 
closed formulae for $M_m(\mathcal{R}_{k})$ are possible:
 
\begin{theorem} Let $m \geq 0$, $k \geq 1$ and $\textbf{\textit{Q}}=\textbf{\textit{P}}_{\overline{M}M} \cdot\textbf{\textit{P}}_{M\overline{M}}$.
\begin{align}
  \text{a) }  & \text{If }  \textbf{\textit{P}} \text{ is } M- \text{ and } \overline{M}-commutable, \text{ then } \nonumber \\
		& M_m(\mathcal{R}_{k}) = k^m \cdot \textbf{\textit{P}}_{M}^{k} 	
		+ \sum_{r=1}^{k} \binom{k}{r} \cdot \textbf{\textit{P}}_{M}^{k-r} 		
		\cdot \textbf{\textit{P}}_{M\overline{M}} \cdot \textbf{\textit{Q}}^{r-1} \cdot  
			\nonumber \\
		& \qquad  \cdot \sum_{j=0}^{m} \binom{j+r-1}{j} \cdot b_{m,j,k+r}
		 \cdot \textbf{\textit{P}}_{\overline{M}}^j 
		 \cdot (\textbf{\textit{I}}-\textbf{\textit{P}}_{\overline{M}})^{-j-r} \cdot \textbf{\textit{P}}_{\overline{M}M}. \label{n3_18} \\
	\text{b) }  & \text{If }  |M|=1 \text{ with } \textbf{\textit{P}}_{M}=(1-p) 
							\text{ and } 
							\textbf{\textit{P}}_{\overline{M}} \cdot \textbf{\textit{e}}_{\overline{M}} =
							s_{\overline{M}} \cdot \textbf{\textit{e}}_{\overline{M}}, \nonumber \\
							& M_m(\mathcal{R}_{k}) = \sum_{r=0}^{k} \binom{k}{r} \cdot 
							\left( \frac{p}{1-p}\right)^r \cdot (1-p)^k \cdot \nonumber \\
			        & \qquad \qquad  \qquad \cdot \sum_{j=0}^{m} \binom{j+r-1}{j} \cdot b_{m,j,k+r} 						 \cdot \left( \frac{s_{\overline{M}}}{1-s_{\overline{M}}}\right) ^j. \label{n3_19}\\
	\text{c) }	& \text{If }  |\overline{M}|=1 \text{ with } 
							\textbf{\textit{P}}_{\overline{M}}=(0) \text{ and } 
							\textbf{\textit{P}}_{M} \cdot \textbf{\textit{e}}_{M} =
							s_{M} \cdot \textbf{\textit{e}}_{M}, \nonumber \\	
							& M_m(\overline{\mathcal{R}}_{k}) = \sum_{j=0}^{m} \binom{j+k-1}{j} \cdot 
							b_{m,j,2 \cdot k} \cdot \left( \frac{s_M}{1-s_M} \right) ^j. \label{n3_20}
\end{align}	
\end{theorem}
\textit{Proof:} a) Induction on $k$. The case $k=1$ is contained in (\ref{n3_15}). If (\ref{n3_18}) holds as induction hypotheses for $k$, then we replace in 
\begin{align*}
	& M_m(\mathcal{R}_{k+1})  
	=\sum_{j=0}^{m} \binom{m}{j} \cdot M_j(\mathcal{R}_{k}) \cdot M_{m-j}(\mathcal{R}_{1})
\end{align*}
(cf. (\ref{n3_16})) $M_j(\mathcal{R}_{k})$ and $M_{m-j}(\mathcal{R}_{1})$ by the formulae 
(\ref{n3_16}) and (\ref{n3_18}), and we get the terms $T_1, T_2, T_3$ and $T_4$ with  
\begin{align*}
T_1=\sum_{j=0}^{m} \binom{m}{j} \cdot k^j \cdot \textbf{\textit{P}}_{M}^k \cdot 
\textbf{\textit{P}}_{M} = (k+1)^m \cdot \textbf{\textit{P}}_{M}^{k+1}, 
\end{align*}
\begin{align*}
&T_2=\\
&= \sum_{j=0}^{m} \binom{m}{j} \cdot k^j \cdot \textbf{\textit{P}}_{M}^k \cdot \textbf{\textit{P}}_{M\overline{M}} \cdot  
	\sum_{i=0}^{m-j} b_{m-j,i,2} \cdot \textbf{\textit{P}}_{\overline{M}}^{i} \cdot (\textbf{\textit{I}}-\textbf{\textit{P}}_{\overline{M}})^{-i-1} \cdot \textbf{\textit{P}}_{\overline{M}M}= \\
	&=\sum_{i=0}^{m} b_{m,i,k+2} \cdot \textbf{\textit{P}}_{M}^k \cdot \textbf{\textit{P}}_{M\overline{M}} \cdot  \textbf{\textit{P}}_{\overline{M}}^{i} \cdot 
	(\textbf{\textit{I}}-\textbf{\textit{P}}_{\overline{M}})^{-i-1} \cdot 
	\textbf{\textit{P}}_{\overline{M}M},
\end{align*}
due to (\ref{a21}),
\begin{align*}
T_3= &\sum_{r=1}^{k} \binom{k}{r} \cdot 
	\sum_{i=0}^{m}	\binom{i+r-1}{i} \cdot b_{m,i,k+1+r}	\cdot \\
	& \qquad \qquad \qquad  \cdot \textbf{\textit{P}}_{M}^{k-r+1} \cdot 
	\textbf{\textit{P}}_{M\overline{M}} \cdot 
	\textbf{\textit{Q}}^{r-1} \cdot \textbf{\textit{P}}_{\overline{M}}^i 
		 \cdot (\textbf{\textit{I}}-\textbf{\textit{P}}_{\overline{M}})^{-i-r} \cdot 
		\textbf{\textit{P}}_{\overline{M}M}, 
\end{align*}
due to $M$-commutability and  (\ref{a21}) and
\begin{align*}
T_4= &\sum_{j=0}^{m} \binom{m}{j} \cdot ( \sum_{r=1}^{k} \binom{k}{r} 
	\cdot \textbf{\textit{P}}_{M}^{k-r} 		
	\cdot \textbf{\textit{P}}_{M\overline{M}} \cdot \textbf{\textit{Q}}^{r-1} \cdot  \\
		&\cdot \sum_{i=0}^{j} \binom{i+r-1}{i} \cdot b_{j,i,k+r}
		 \cdot \textbf{\textit{P}}_{\overline{M}}^i 
		 \cdot (\textbf{\textit{I}}-\textbf{\textit{P}}_{\overline{M}})^{-i-r} \cdot \textbf{\textit{P}}_{\overline{M}M}) \cdot \\
		 & \qquad \qquad \qquad \cdot (\textbf{\textit{P}}_{M\overline{M}} \cdot  
	\sum_{i=0}^{m-j} b_{m-j,i,2} \cdot \textbf{\textit{P}}_{\overline{M}}^{i} \cdot (\textbf{\textit{I}}-\textbf{\textit{P}}_{\overline{M}})^{-i-1} \cdot \textbf{\textit{P}}_{\overline{M}M})= \\
	=&\sum_{r=1}^{k} \binom{k}{r} \cdot \sum_{i=0}^{m} b_{m,i,k+2+r} \cdot
	\textbf{\textit{P}}_{M}^{k-r} 		
	\cdot \textbf{\textit{P}}_{M\overline{M}} \cdot \textbf{\textit{Q}}^{r-1} \cdot  \\
	& \qquad \qquad \cdot \sum_{s=0}^{i} \binom{s+r-1}{s} \cdot \textbf{\textit{P}}_{\overline{M}}^{s} \cdot (\textbf{\textit{I}}-\textbf{\textit{P}}_{\overline{M}})^{-s-r} \cdot \textbf{\textit{P}}_{\overline{M}M} \cdot \textbf{\textit{P}}_{M\overline{M}} \cdot \\
	&\qquad \qquad \qquad \qquad \textbf{\textit{P}}_{\overline{M}}^{i-s} \cdot (\textbf{\textit{I}}-\textbf{\textit{P}}_{\overline{M}})^{s-i-1} \cdot \textbf{\textit{P}}_{\overline{M}M} = \\
	=&\sum_{r=1}^{k} \binom{k}{r} \cdot \sum_{i=0}^{m} b_{m,i,k+2+r} \cdot
	\textbf{\textit{P}}_{M}^{k-r} 		
	\cdot \textbf{\textit{P}}_{M\overline{M}} \cdot \textbf{\textit{Q}}^{r} \cdot  \\
	& \qquad \qquad \cdot \binom{i+r}{i} \cdot \textbf{\textit{P}}_{\overline{M}}^{i} \cdot (\textbf{\textit{I}}-\textbf{\textit{P}}_{\overline{M}})^{-i-r-1} \cdot \textbf{\textit{P}}_{\overline{M}M},  
\end{align*}
due to $\overline{M}$-commutability and (\ref{a18}). Finally
\begin{align*}
&M_m(\mathcal{R}_{k+1}) = T_1 + T_2 + T_3 + T_4 = \\
& =(k+1)^m \cdot \textbf{\textit{P}}_{M}^{k+1}  
 +\sum_{r=1}^{k+1} \binom{k+1}{r} \cdot \textbf{\textit{P}}_{M}^{k+1-r} 
	  \sum_{i=0}^{m}	\binom{i+r-1}{i} \cdot b_{m,i,k+1+r}	\cdot \\
		& \qquad \qquad \cdot 
	  \textbf{\textit{P}}_{M\overline{M}} \cdot 
	  \textbf{\textit{Q}}^{r-1} \cdot \textbf{\textit{P}}_{\overline{M}}^i 
		\cdot (\textbf{\textit{I}}-\textbf{\textit{P}}_{\overline{M}})^{-i-r} \cdot 
		\textbf{\textit{P}}_{\overline{M}M}.
\end{align*}

		b) follows from a), because \textbf{\textit{P}} is $M-$commutable and the $\overline{M}-$ commutability was only used to show in the computation of $T_4$ that 
		
		\begin{align*}
		&\textbf{\textit{P}}_{M\overline{M}} \cdot \textbf{\textit{Q}}^{r-1} \cdot	
		\textbf{\textit{P}}_{\overline{M}}^s 
		\cdot (\textbf{\textit{I}}-\textbf{\textit{P}}_{\overline{M}})^{-s-r} \cdot 
		\textbf{\textit{P}}_{\overline{M}M} \cdot  
		\textbf{\textit{P}}_{M\overline{M}} \cdot \\
		& \qquad \qquad \qquad \qquad \textbf{\textit{P}}_{\overline{M}}^{i-s} \cdot  
		(\textbf{\textit{I}}-\textbf{\textit{P}}_{\overline{M}})^{s-i-1} \cdot 
		\textbf{\textit{P}}_{\overline{M}M}= \\
		&\textbf{\textit{P}}_{M\overline{M}} \cdot \textbf{\textit{Q}}^{r} \cdot	
		\textbf{\textit{P}}_{\overline{M}}^i 
		\cdot (\textbf{\textit{I}}-\textbf{\textit{P}}_{\overline{M}})^{-i-r-1} \cdot 
		\textbf{\textit{P}}_{\overline{M}M}. 
		\end{align*}
		This equality holds because of
			\begin{align}
				&\textbf{\textit{P}}_{M\overline{M}} \cdot \textbf{\textit{P}}_{\overline{M}M} =
				p \cdot (1-s_{\overline{M}}), \label{n3_21} \\
				&\textbf{\textit{P}}_{\overline{M}}^n \cdot \textbf{\textit{P}}_{\overline{M}M} =
				s_{\overline{M}}^n \cdot \textbf{\textit{P}}_{\overline{M}M} 
				\label{n3_22} 
				\text{ and} \\
				&(\textbf{\textit{I}}-\textbf{\textit{P}}_{\overline{M}})^{-n} \cdot 
		    \textbf{\textit{P}}_{\overline{M}M} =
				(1-s_{\overline{M}})^{-n} \cdot \textbf{\textit{P}}_{\overline{M}M} 
				\text{  for } n \ge 0 .\label{n3_23} 
			\end{align}
			and equals $\left( \frac{p}{1-p}\right)^r \cdot (1-p)^k \cdot 
			\left( \frac{s_{\overline{M}}}{1-s_{\overline{M}}}\right) ^j$. 
		c) follows from (\ref{n3_18}) by exchanging $M$ and $\overline{M}$ similarly to (\ref{n3_19}). $\qed$

		For $|M|=1$, item b) of the preceding theorem describes the moments of the $k-$th recurrence to a single state of a DTMC and item c) describes the moments of the time to the $k-$th renewal $\overline{\mathcal{R}}_{k}$.  
		
		The next theorem generalizes the negative binomial distribution with the special cases $|\overline{M}|=1$ (phase type distribution) and $|M|=|\overline{M}|=1$ (alternating negative binomial distribution).  
\begin{theorem} Let $m \geq 0$, $k \geq 1$.
\begin{align}
 \text{a) }    &\text{If }  \textbf{\textit{P}} \text{ is } M- \text{ and } \overline{M}-commutable \text{ and } \textbf{\textit{Q}}=\textbf{\textit{P}}_{\overline{M}M} \cdot \textbf{\textit{P}}_{M\overline{M}}, \text{ then } \nonumber \\
 		 &M_m(\mathcal{N}_{k}) =\sum_{r=0}^{k-1} \binom{k-1}{r} \sum_{j=0}^{m} 
		b_{m,j,k+r} \cdot \label{n3_25} \\
		   &\qquad \qquad \cdot \binom{j+r}{j}  \cdot \textbf{\textit{P}}_{M}^j \cdot 
		(\textbf{\textit{I}}-\textbf{\textit{P}}_{M})^{-(j+r+1)} \cdot 
		\textbf{\textit{P}}_{M\overline{M}} \cdot \textbf{\textit{P}}_{\overline{M}}^{k-1-r} \cdot
		\textbf{\textit{Q}}^{r}. \nonumber\\
   \text{b) }  &  \text{If }	|\overline{M}|=1 \text{ with } 
	\textbf{\textit{P}}_{\overline{M}}=(q) \text{ and } 
	\textbf{\textit{P}}_{M} \cdot \textbf{\textit{e}}_{M} =s_M \cdot \textbf{\textit{e}}_{M}, 
	\text{ then } \nonumber \\
 	&  M_m(\mathcal{N}_{k}) =\sum_{r=0}^{k-1} \binom{k-1}{r} \cdot \left(\frac{1-q}{q}\right)^r 
	\cdot q^{k-1} \cdot \sum_{j=0}^{m} 
		b_{m,j,k+r} \cdot \label{n3_26} \\
		& \qquad \qquad \cdot \binom{j+r}{j} \cdot \left(\frac{s_{M}}{1-s_{M}}\right)^j \cdot
		\textbf{\textit{e}}_{M}. \nonumber\\
 	\text{c) }  &  \text{If }	|\overline{M}|=|M|=1 \text{ with } 
	\textbf{\textit{P}}_{\overline{M}}=(q) \text{ and } \textbf{\textit{P}}_{M}=(1-p), 
	\text{ then } \nonumber \\
  		& M_m(\mathcal{N}_{k}) =\sum_{r=0}^{k-1} \binom{k-1}{r} \left(\frac{1-q}{q}\right)^r 
	\cdot q^{k-1} \sum_{j=0}^{m} 
		b_{m,j,k+r} \binom{j+r}{j} \cdot \left(\frac{1-p}{p}\right)^j. \label{n3_27}  \\
 	\text{d) }  &  \text{If in c) } p=q, \text{ the moments of the negative binomial distribution are}	\nonumber \\
  	& M_m(\mathcal{N}_{k}) =\sum_{j=0}^{m} \binom{j+k-1}{k-1} \cdot b_{m,j,k} \cdot \left(\frac{1-p}{p}\right)^j. \label{n3_28}  
\end{align} 
\end{theorem} 
\textit{Proof:} a) Due to (\ref{n3_17}), (\ref{n3_14}), (\ref{n3_18}) , (\ref{a21}) and the commutability, with $\overline{\textbf{\textit{Q}}}=\textbf{\textit{P}}_{M\overline{M}} \cdot \textbf{\textit{P}}_{\overline{M}M}$
\begin{align*}
&M_m(\mathcal{N}_k) = 
	     \sum_{j=0}^{m} \binom{m}{j}\cdot M_{m-j}(\mathcal{N}_{1}) 
			 \cdot M_{j}(\mathcal{\overline{R}}_{k-1})= \\
			&=\sum_{j=0}^{m} \binom{m}{j} \cdot 
			(\sum_{r=0}^{m-j} b_{m-j,r,1} \cdot \textbf{\textit{P}}_M^{r} \cdot (\textbf{\textit{I}}-\textbf{\textit{P}}_M)^{-r-1} \cdot \textbf{\textit{P}}_{M\overline{M}}) \cdot \\
			&\qquad \cdot ((k-1)^j \cdot \textbf{\textit{P}}_{\overline{M}}^{k-1}+ \sum_{r=1}^{k-1} 
			\binom{k-1}{r} \cdot \textbf{\textit{P}}_{\overline{M}}^{k-1-r} 		
			\cdot \textbf{\textit{P}}_{\overline{M}M} \cdot \overline{\textbf{\textit{Q}}}^{r-1} \cdot  \\
		  &\qquad \cdot \sum_{i=0}^{j} \binom{i+r-1}{i} \cdot b_{j,i,k-1+r}
		 \cdot \textbf{\textit{P}}_{M}^i 
		 \cdot (\textbf{\textit{I}}-\textbf{\textit{P}}_{M})^{-i-r} \cdot \textbf{\textit{P}}_{M\overline{M}})= \\
		&=\sum_{r=0}^{m} \sum_{j=0}^{m-r} \binom{m}{j} \cdot b_{m-j,r,1} \cdot (k-1)^j \cdot
		\textbf{\textit{P}}_M^{r} \cdot (\textbf{\textit{I}}-\textbf{\textit{P}}_M)^{-r} \cdot
		\textbf{\textit{P}}_{M\overline{M}} \cdot \textbf{\textit{P}}_{\overline{M}}^{k-1} + \\
		&\qquad +\sum_{l=0}^{m} \sum_{r=0}^{l} \textbf{\textit{P}}_M^{r} \cdot 
		(\textbf{\textit{I}}-\textbf{\textit{P}}_M)^{-r} \cdot 
		\textbf{\textit{P}}_{M\overline{M}} \cdot \sum_{s=0}^{k-2} \binom{k-1}{s} \cdot \\ 
		&\qquad \cdot \sum_{j=0}^{m} \binom{m}{j} \cdot b_{m-j,r,1} \cdot 
		b_{j,l-r,2 \cdot (k-1)-s} \cdot  \binom{l-r+k-2-s}{k-2-s} \cdot 
		\textbf{\textit{P}}_{\overline{M}}^s \cdot \\
		&\qquad \cdot 
		\textbf{\textit{P}}_{\overline{M}M} \cdot
		\overline{\textbf{\textit{Q}}}^{k-2-s} \cdot 
		(\textbf{\textit{I}}-\textbf{\textit{P}}_M)^{-(k-1-s)}\cdot\textbf{\textit{P}}_M^{l-r}\cdot
		(\textbf{\textit{I}}-\textbf{\textit{P}}_M)^{r-l}\cdot \textbf{\textit{P}}_{M\overline{M}}=\\
		&=\sum_{r=0}^{m} b_{m,r,k} \cdot
		\textbf{\textit{P}}_M^{r} \cdot (\textbf{\textit{I}}-\textbf{\textit{P}}_M)^{-r-1} \cdot
		\textbf{\textit{P}}_{M\overline{M}} \cdot \textbf{\textit{P}}_{\overline{M}}^{k-1} + \\
		&\qquad +\sum_{l=0}^{m} \sum_{s=0}^{k-2} \binom{k-1}{s} \cdot \binom{l+k-s-1}{k-s-1} \cdot
		b_{m,l,2 \cdot k -1-s} \cdot \\
		&\qquad \cdot \textbf{\textit{P}}_M^{l} \cdot 
		(\textbf{\textit{I}}-\textbf{\textit{P}}_M)^{-(l+k-s)} \cdot
		\textbf{\textit{P}}_{M\overline{M}} \cdot \textbf{\textit{P}}_{\overline{M}}^s 
		\cdot \textbf{\textit{Q}}^{k-1-s} = \\
		& =\sum_{l=0}^{m} \sum_{s=0}^{k-1} \binom{k-1}{s} \cdot \binom{l+k-s-1}{k-s-1} \cdot
		b_{m,l,2 \cdot k -1-s} \cdot \\
		&\qquad \cdot \textbf{\textit{P}}_M^{l} \cdot 
		(\textbf{\textit{I}}-\textbf{\textit{P}}_M)^{-(l+k-s)} \cdot
		\textbf{\textit{P}}_{M\overline{M}} \cdot \textbf{\textit{P}}_{\overline{M}}^s 
		\cdot \textbf{\textit{Q}}^{k-1-s},
\end{align*}
and the index transformation $r=k-s-1$ leads to (\ref{n3_25}).
b) follows from a) due to (\ref{n3_21}), (\ref{n3_22}) and (\ref{n3_23}). c) follows from b) due to $s_M=1-p$. d) Induction on $k$. For $k=1$, (\ref{n3_28}) is (\ref{n3_14}). For the step  from $k$ to $k+1$, we remark that under the assumption $p=q$, the distributions of $\mathcal{N}_{1}$ and $\overline{\mathcal{R}}_{1}$ are equal and therefore also the distributions of $\mathcal{N}_{k}$ and $\overline{\mathcal{R}}_{k}$. 
So, $\mathcal{N}_{k+1}=\mathcal{N}_{k}+\mathcal{N}_{1}$ and due to (\ref{a17})
\begin{align*}
	&M_m(\mathcal{N}_{k+1})=
	\sum_{r=0}^{m} \binom{m}{r} \cdot M_{m-r}(\mathcal{N}_{k}) \cdot M_r(\mathcal{N}_{1})=\\
	&\sum_{r=0}^{m} \binom{m}{r} \cdot
	\sum_{s=0}^{m-r} \binom{s+k-1}{k-1 }\cdot b_{m-r,s,k} \cdot \left( \frac{1-p}{p}\right)^s
	 \cdot \sum_{s=0}^{r}  b_{r,s,1} \cdot \left( \frac{1-p}{p}\right)^s=\\
	&\sum_{r=0}^{m} \binom{m}{r} \cdot \sum_{l=0}^{m} \sum_{s=0}^{l} 
	b_{m-r,s,k} \cdot b_{r,l-s,1} \cdot \left( \frac{1-p}{p}\right)^l \cdot \binom{s+k-1}{k-1}=\\
	&\sum_{l=0}^{m} \left( \frac{1-p}{p}\right)^l b_{m,l,k+1}  \cdot
	\sum_{s=0}^{l} \binom{s+k-1}{k-1}= 
	\sum_{l=0}^{m}  b_{m,l,k+1} \cdot \binom{l+k}{k} \cdot \left( \frac{1-p}{p}\right)^l . \\
	 \square{}
\end{align*}
	
Of course, also the binomial distribution can be generalized like the negative binomial distribution to an "`alternating binomial distribution"', and its moments can be computed in an analogous way.

Another group of applications deals with the central moments 
	\[C_m(X) = M_1((X-M_1(X))^m),\] 
if $X$ is a random variable on $\{0,1,2, \ldots\}$. If $F_m(X)= M_1(X\cdot(X-1)\cdot \ldots \cdot (X-m+1))$ is the $m-$th factorial moment, then by means of (\ref{a7a}) 
\begin{align}
 &M_m(X)=\sum_{j=0}^{m} b_{m,j,0} \cdot \frac{1}{j!} \cdot F_j(X) \text{ and}  \label{n_29}\\
 &C_m(X)=\sum_{j=0}^{m} \binom{m}{j} \cdot(-M_1(X))^{m-j} \cdot M_j (X).\label{n_30}
\end{align} 
The next theorem connects $C_m(X)$ with $F_m(X)$.
\begin{theorem} If $X$ is a random variable on $\{0,1,2, \ldots\}$ with existing moments of any order and $m \ge 0$,
\begin{align}
C_m(X)=\sum_{j=0}^{m} b_{m,j,-M_1(X)} \cdot \frac{1}{j!} \cdot F_j(X).\label{n_3_31}
\end{align}
\end{theorem}
\textit{Proof.} From (\ref{n30}) follows due to (\ref{a20}) and (\ref{n_29}) 
\begin{align*}
C_m(X)=&\sum_{j=0}^{m} \binom{m}{j} \cdot (-M_1(X))^{m-j} \cdot M_j (X)= \\
&\sum_{j=0}^{m} \binom{m}{j} \cdot (-M_1(X))^{m-j} \cdot \sum_{i=0}^{j} b_{j,i,0} \cdot 
\frac{1}{i!} \cdot F_i(X)=\\
&\sum_{i=0}^{m}\frac{1}{i!} \cdot F_i(X) \cdot 
\sum_{j=0}^{m} \binom{m}{j} \cdot (-M_1(X))^{m-j} \cdot b_{j,i,0} =\\
&\sum_{i=0}^{m} b_{m,i,-M_1(X)} \cdot \frac{1}{i!} \cdot F_i(X). \qed 
\end{align*}
Similar transforms from (ordinary) moments to central moments are possible for the binomial, the Poisson, the negative binomial, the alternating binomial, the discrete uniform distribution (equal probabilities for the values $0,1,\ldots N-1$ and the phase type distribution. We recall for the random variables $X_{n,p}\sim B(n,p)$, $X_{\lambda}\sim Poi(\lambda)$, $X_{p,k}\sim NB(p,k)$, $X_{p,q,k}\sim ANB(p,q,k)$, $X_{N}\sim Hom(N)$, 
$X_{\textbf{\textit{a}},\textbf{\textit{A}}} \sim PH_D(\textbf{\textit{a}},\textbf{\textit{A}})$ with initial vector $\textbf{\textit{a}}$ and matrix $\textbf{\textit{A}}$ and the recurrence time $\mathcal{R}_1$ of state 1 of the DTMC defined by the matrix $\textbf{\textit{P}}$ with $M=\{1\}$, then 
\begin{equation}
M_m(X_{n,p})=	\sum_{j=0}^{n} S_{m,j} \cdot (n)_{j} \cdot p^j =
	\sum_{j=0}^{n} b_{m,j,0} \cdot \binom{n}{j} \cdot p^j \label{a48}
\end{equation}
(cf. e.g. \cite{Binomi}),
\begin{equation}
M_m(X_{\lambda})=\sum_{n=0}^{\infty} n^m \cdot P(X_{\lambda}=n)= 
	\sum_{j=0}^{m} S_{m,j} \cdot \lambda ^j =	
	\sum_{j=0}^{m} b_{m,j,0} \cdot \frac{\lambda ^j}{j!} \label{a47}
\end{equation}

(cf. e.g. \cite{Stirling}), the formulas (\ref{n3_28}) and (\ref{n3_27}) for $M_m(X_{p,k})$, resp. $M_m(X_{p,q,k})$, 
\begin{equation} 
M_m(X_N)= \frac{1}{N} \cdot \sum_{j=0}^{m} b_{m,j,0} \cdot \binom{N}{r+1}
\end{equation}
and the formula (\ref{n3_15}) for $PH_D(\textit{\textbf{a}},\textit{\textbf{A}})$ and the recurrence time $\mathcal{R}_1$.
The next theorem shows that only the third parameter of the \textit{MSNs} changes, if central moments are computed instead of die ordinary ones.
\begin{theorem} With $m \ge 0$, $X_{n,p}$, $X_{\lambda}$, $X_{p,k}$, $X_{p,q,k}$, $X_{N}$, $X_{\textbf{\textit{a}},\textbf{\textit{A}}}$ and $\mathcal{R}_1$ as defined above 
\begin{align}
	&\text{a) } C_m(X_{n,p})=\sum_{j=0}^{n} b_{m,j,-n \cdot p} \cdot \binom{n}{j} \cdot p^j ,\\
	&\text{b) } C_m(X_{\lambda})=\sum_{j=0}^{m} b_{m,j,-\lambda} \cdot \frac{\lambda ^j}{j!} ,\\
	&\text{c) } C_m(X_{p,k}) =\sum_{j=0}^{m} \binom{j+k-1}{k-1} \cdot b_{m,j,k 
	\cdot (1-\frac{1}{p})} \cdot 	\left(\frac{1-p}{p}\right)^j,\\
	&\text{d) } C_m(X_{p,q,k}) =\sum_{r=0}^{k-1} \binom{k-1}{r} \cdot \left(\frac{1-q}{q}\right)^r 
	\cdot q^{k-1} \cdot \\
	& \qquad \qquad \qquad \qquad \cdot \sum_{j=0}^{m} b_{m,j,k+r-M_1(X_{p,q,k})} \cdot
	\binom{j+r}{j} \cdot 
	\left(\frac{1-p}{p}\right)^j \nonumber \\ 
	&\qquad \qquad \qquad \qquad \text{ with } M_1(X_{p,q,k})= \frac{(k-1)\cdot (p-q) +k}{p},
	 \nonumber \\
	&\text{e) } C_m(X_N)= \frac{1}{N} \cdot \sum_{j=0}^{m} b_{m,j,-\frac{N-1}{2}} 
	\cdot \binom{N}{j+1}, \\
	&\text{f) } C_m(X_{\textbf{\textit{a}},\textbf{\textit{A}}})=(1-\textbf{\textit{a}} \cdot 
	\textbf{\textit{e}}_{dim(\textbf{\textit{A}})}) \cdot
	\left(1-M_1(X_{\textbf{\textit{a}},\textbf{\textit{A}}})\right)^m + \\
	&\qquad \qquad + \textbf{\textit{a}} \cdot \sum_{j=0}^{m} b_{m,j,2-M_1(X_{\textbf{\textit{a}},\textbf{\textit{A}}})} \cdot \textbf{\textit{A}}^j \cdot 
	(\textbf{\textit{I}}-\textbf{\textit{A}})^{-j} \cdot 
	\textbf{\textit{e}}_{dim(\textbf{\textit{A}})} \nonumber \\
	&\qquad \qquad \text{ with } M_1(X_{\textbf{\textit{a}},\textbf{\textit{A}}})= 
	1 + \textbf{\textit{a}} \cdot (\textbf{\textit{I}}-\textbf{\textit{A}})^{-1} \cdot \textbf{\textit{e}}_{dim(\textbf{\textit{A}})} ,\nonumber \\
	&\text{g) } C_m(\mathcal{R}_1)= 
	\textbf{\textit{P}}_{M} \cdot \left(1-M_1(\mathcal{R}_1\right)^m +  \\
	& \qquad \qquad + \textbf{\textit{P}}_{M\overline{M}} \cdot  
	\sum_{j=0}^{m} b_{m,j,2-M_1(\mathcal{R}_1)} \cdot \textbf{\textit{P}}_{\overline{M}}^{j} \cdot (\textbf{\textit{I}}-\textbf{\textit{P}}_{\overline{M}})^{-j-1} \cdot 
	\textbf{\textit{P}}_{\overline{M}M}  \nonumber \\
	& \qquad \qquad \text{ with } M_1(\mathcal{R}_1) = 1+ \textbf{\textit{P}}_{M\overline{M}}\cdot 
	(\textbf{\textit{I}}-\textbf{\textit{P}}_{\overline{M}})^{-1} \cdot \textbf{\textit{e}}_{\overline{M}}. \nonumber
\end{align}
\end{theorem}  
The proofs follow from (\ref{n_30}) with (\ref{a20}) similarly to (\ref{n_3_31}).
 
Some questions must remain open here, which should be answered in further investigations:
\begin{itemize}
\item What is a combinatorial interpretation of \textit{MSNs} and how can they be applied, if $k$ is not an integer greater or equal zero?
\item What is a combinatorial interpretation of \textit{MSN1s} and how can they be applied?
\item What are the upper and lower bounds for \textit{MSNs} and \textit{MSN1s}?
\item What are the relationships between the \textit{MSN1s} and to the \textit{MSNs}?
\item Are there other applications of the \textit{MSNs}?
\end{itemize}

\renewcommand{\refname}{References}
\setlength{\bibitemsep}{0.0\baselineskip}

\end{document}